\documentclass[12pt]{article}
\usepackage{amssymb,amsthm}
\small\normalsize

\begin{document}

\newtheorem{theorem}{Theorem}
\newtheorem{corollary}[theorem]{Corollary}
\newtheorem{lemma}[theorem]{Lemma}

\def\ex{{\bf ex}}

\pagestyle{myheadings}
\markright{{\small{\sc Z. F\"uredi:} Graphs not containing a cube}}

\title{ \phantom{~}\vskip-1.8cm
\bf On a theorem of Erd\H os and Simonovits on graphs
    not containing the cube}
\author{{\bf Zolt\'an F\"uredi}}
\date{\small ${}$
R\'enyi Institute of Mathematics of the Hungarian Academy of Sciences
  \\ Budapest, P. O. Box 127, Hungary-1364. \\ \medskip
E-mail: {\tt z-furedi@illinois.edu, furedi.zoltan@renyi.mta.hu}
\footnotetext{\noindent
Research supported in part by the Hungarian National Science Foundation
OTKA 104343,
and by the European Research Council Advanced Investigators Grant 267195.\\
MSC-class: 05C35, 05D99
}}
\maketitle

\parskip=3pt
\begin{abstract}
The cube $Q$ is the usual $8$-vertex graph with $12$ edges.
Here we give a new proof for a theorem
 of Erd\H os and Simonovits concerning the Tur\'an number of the cube.
Namely, it is shown that $e(G) \leq  n^{8/5}+(2n)^{3/2}$ holds for any $n$-vertex
 cube-free graph $G$.

Our aim is to give a self-contained exposition. We also
 point out the best known results and supply bipartite versions.
\end{abstract}

\normalsize
\section{History of Tur\'an type problems}

As usual, we write $|G|$, $e(G)$, $\deg_G(x)$ for the number of vertices, number of edges, and the degree
 of a vertex $x$ of a graph $G$.
Denote by $N_G(x)$ (or just $N(x)$) the neighborhood of $x$, note that $x\notin N(x)$.
Let $K_n$ and $K_{a,b}$ denote the complete graph on $n$ vertices and the complete
 bipartite graph with bipartition classes of sizes $a$ and $b$.
$K(A,B)$ denotes the complete bipartite graph with partite sets $A$ and $B$ ($A\cap B=\emptyset$).

A graph $G$ not containing $H$ as a (not necessarily induced) subgraph is called \emph{$H$-free}.
Let us denote by $\ex(n,H)$ the \emph{Tur\'an number} for $H$,
 i.e.~the maximum number of edges of an $H$-free graph on $n$ vertices.
More generally, let $\ex(G,H)$ be the maximum number of edges in an $H$-free subgraph of $G$.
Then $\ex(n,H)=\ex(K_n,H)$.
We also use the notation $\ex(a,b,H)$ for $\ex(K_{a,b}, H)$ and call it the bipartite version of Tur\'an number.
Also, if $F \subset H$ then $\ex(n,F) \leq \ex(n,H)$.

Tur\'an~\cite{TuranML} determined  $\ex(n,K_{p+1})$.
The extremal graph is the
 almost equipartite complete graph of $p$ classes.
He also proposed the general question, $\ex(n,H)$, in particular the determination of
 the Tur\'an number of the graphs obtained from the platonic polyhedrons,
 the cube $Q=Q_8$ (it is an 8-vertex 3-regular graph), the octahedron $O_6$
  (six vertices, 12 edges), the icosahedron $I_{12}$ (12 vertices, 5-regular) and for the
  dodecahedron $D_{20}$ (20 vertices, 30 edges).
Erd\H{o}s and Simonovits~\cite{ErdSimOcta} gave an implicit formula for
  $\ex(n, O_6)$ (they reduced it to $\ex(n,C_4)$) and Simonovits solved
  exactly  $\ex(n, D_{20})$ in~\cite{SimSymm} and $\ex(n , I_{12})$ in~\cite{SimIco}
  (for $n> n_0$).

In fact, Tur\'an's real aim was not only these particular graphs but to discover
 a general theory.
His questions, and the answers above, indeed lead to an asymptotic
 (the Erd\H os-Simonovits theorem~\cite{ErdSimLim}) and to the Simonovits stability
 theorem concerning the extremal graphs~\cite{SimTihany} in the case when
 the sample graph $H$ has chromatic number at least three.
For a survey and explanation see  Simonovits~\cite{SimFra} or the monograph of Bollob\'as~\cite{Bolkonyv}.

However, the bipartite case is different, see the recent survey~\cite{FureSimoErdos100}.
Even the extremal problem of the cube graph $Q$,
 which was one of Tur\'an's  \cite{TuranColloq} originally posed problems, is still unsolved.
Our aim here is to give a gentle introduction to this topic.
We survey the results and methods concerning $\ex(n,Q)$,  give new or at least
 streamlined proofs.
We only use basic ideas of multilinear optimization (Lagrangian, convexity, etc.) and in most cases
 just high school algebra.
We also consider the case of bipartite host graph, i.e., $\ex(a,b,Q)$.

\section{Walks}

Let $W_3=W_3(G)$ denote the number of walks in $G$ of length 3, i.e.,
 the number of sequences of the form $x_0x_1x_2x_3$ where $x_{i-1}x_{i}$ is
 an edge of $G$ (for $i=1,2,3$).
Note that, e.g., $xyxy$ is a walk (if $xy\in E(G)$) and it differs from
 $yxyx$.
A $d$-regular graph has exactly $nd^3$ 3-walks.

\begin{theorem}\label{W3th}
For every $n$-vertex graph $G$ for the number of 3-walks one has
\begin{equation}\label{W3eq1}
 W_3\geq n\left(\frac{1}{n}\sum_{x\in V} \deg(x)^{3/2}\right)^2 .
 \end{equation}
 \end{theorem}

The $r$-order {\bf power mean} of the nonnegative sequence $a_1, \dots, a_m $
 is $M_r({\bf a}):= \left(\frac{1}{m}\sum a_i^r \right)^{1/r}$.
Then for $1\leq r\leq  s\leq \infty$ one has
\begin{equation}\label{cubeeq1}
a_{\rm ave}:=M_1({\bf a})\leq M_r({\bf a})\leq M_s({\bf a})\leq
  M_{\infty}({\bf a}):= \max_i |a_i|.
  \end{equation}
We will frequently use it in the equivalent form
\begin{equation}\label{eq:3}
\sum_{1\leq i \leq m} a_i^r \leq \left(\sum a_i^s\right)^{r/s} m^{1- (r/s)}.
  \end{equation}
This is just a special case of the H\"older inequality, i.e.,
 for any two non-negative vectors ${\bf x}, {\bf y}\in R^{m}$
 and for reals $p,q\geq 1$ with $\frac{1}{p}+\frac{1}{q}=1$ one has
$$
    \sum_i x_iy_i  \leq \left(\sum_i{x_i^p}\right)^{1/p}
  \left(\sum_i{y_i^q}\right)^{1/q}.
  $$
We get (\ref{eq:3}) by substituting here
${\mathbf x}=(a_i^r)_{1\leq i\leq m}$, ${\bf y}=(1,1,\dots, 1)$, $1/p=r/s$ and $1/q=1-(r/s)$.

\medskip
\noindent{\it Proof of Theorem~\ref{W3th}.}\quad
Considering the middle edge of the 3-walks one obtains that
$$
 W_3=\sum_{x\in V}  \sum_{y \in N(x)} \deg(x)\deg(y) .
 $$
Here we have $2e=nd_{\rm ave}$ terms.
We use for this sum the $2e$-dimensional Chauchy-Schwartz inequality
$$
 \left(\sum_i{a_i^2}\right)\left(\sum_i{b_i^2}\right)\geq \left(\sum_i a_ib_i\right)^2
 $$
 valid for any two vectors ${\bf a}, {\bf b}\in R^{m}$.
Our aim is to separate the variables in the products $\deg(x)\deg(y)$
 so we take ${\bf a}= \{\sqrt{\deg(x)\deg(y)}\}_{x\in V,\,y\in N(x)} $ and
 ${\bf b}=\{ \sqrt{1/\deg(x)}\}_{x\in V,\,y\in N(x)}$.
One obtains that
\begin{eqnarray*}
 W_3\, n &=& \left(\sum_{x\in V}  \sum_{y \in N(x)} \deg(x)\deg(y)\right)
      \left(\sum_{x\in V}  \sum_{y \in N(x)} \frac{1}{\deg(x)}\right)
        = \left(\sum_i{a_i^2}\right)\left(\sum_i{b_i^2}\right)\\
  &\geq&\left(\sum_i a_ib_i\right)^2=
  \left(\sum_{x\in V}  \sum_{y \in N(x)}
   \frac{\sqrt{\deg(x)\deg(y)}}{\sqrt{\deg (x)}} \right)^2
 =\left(\sum_{x\in V}  \sum_{y \in N(x)} \sqrt{\deg(y)} \right)^2 \\
 &=&\left(\sum_{y\in V} \deg(y)^{3/2}\right)^2.
 \hfill \qed  \end{eqnarray*}

\medskip
\noindent{\bf Historical remarks.}\quad
One can rewrite Theorem~\ref{W3th} as
\begin{equation}\label{eq:4}
  W_3\geq nM_{3/2}({\bf d})^3.
  \end{equation}
Then the power mean inequality (\ref{cubeeq1}) with $(r,s)=(1,3/2)$ gives that
\begin{equation}\label{eq:5}
  W_3\geq n (d_{\rm ave})^3 =8e^3/n^2.
    \end{equation}
This inequality $W_3\geq n (d_{\rm ave})^3 $ is due to  Mulholland and
 C. A. B. Smith~\cite{MS} and was generalized
 by Attkinson,  Watterson and Moran~\cite{AWM} for  $W_k$ for every $k\geq 3$
 in a form of a matrix inequality.
Then it was further generalized by
 Blakley and Roy~\cite{BlakleyRoy} for all nonnegative symmetric matrices.
As far as the author knows the obvious consequence of their works,
 $W_k\geq n (d_{\rm ave})^k$, was first explicitly stated in a paper
  of Erd\H os and Simonovits~\cite{ESpath}.
For the interested reader we supply a direct proof for (\ref{eq:5}) using only high school algebra
 in the Appendix (Section~\ref{S:App}).

Theorem~\ref{W3th} is not really new.
It is an easy consequence (of a special case) of a result of
 Jagger,  \v{S}\'tov\'{\i}\v{c}ek, and Thomason~\cite{Andrew},
 who while working on a conjecture of Sidorenko~\cite{Sid} showed
 the inequality
   $\sum_x w(x)^{1/2}\geq \sum_x \deg(x)^{3/2}$
 where $w(x)$ is the number of 3-walks whose second vertex is  $x$.

\noindent{\bf The exponent $3/2$ is the best possible.}\quad
Consider a complete bipartite graph $K_{a,b}$, we have $W_3=2a^2b^2$.
Then $W_3/nM_p({\bf d})^3\to 0$
 for any fixed $p> 3/2$ whenever $b/a\to \infty$.

For $K_{a,b}$ we have $W_3=2a^2b^2$, while the right hand side of (\ref{W3eq1}) is
 $a^2b^2\frac{(\sqrt{a}+\sqrt{b})^2}{a+b}$ which is between
 $a^2b^2$ and $2a^2b^2$.
Using this observation one can show the following:
Suppose that $d_1, \dots, d_n$ is the degree sequence of a
 graph $G$.
There is a graph $H$ with degree sequence $d_i'$
 for which $d_i'\geq d_i$ and
 $W_3(H)\leq 4 n(\sum (d_i')^{3/2}/n)^2=4nM_{3/2}({\bf d'})^3$.

\section{3-paths in bipartite graphs}

Let $P_3=P_3(G)$ denote the number of 3-paths of $G$.
We have $2P_3\leq W_3$.
Using the method of the previous Section we show the following
 lower bound for $P_3$.

\begin{theorem}\label{P3th}
Let $G(A,B)$ be a bipartite graph with $e$ edges and with color classes
 $A$ and $B$, $|A|=a$, $|B|=b$.
Suppose that every vertex has degree at least 2.
Then for the number of 3-paths one has
\begin{equation}\label{eq:P3}
 P_3\geq \frac{e(e-a)(e-b)}{ab}.
 \end{equation}
 \end{theorem}

\noindent{\it Proof.}\quad
Considering the middle edge of the 3-paths one obtains that
$$
 P_3=\sum_{x\in A}  \sum_{y \in N(x)} (\deg(x)-1)(\deg(y)-1) .
 $$
Here we have $e$ terms.
One obtains that
\begin{eqnarray*}
 a\,P_3 &=& a \sum_{x\in A}  \sum_{y \in N(x)} (\deg(x)-1)(\deg(y)-1) \\
 &=& a\sum_{x\in A}  \sum_{y \in N(x)} -(\deg(y)-1) +
     a\sum_{x\in A}  \sum_{y \in N(x)} \deg(x)(\deg(y)-1) \\
=& -a&\!\!\!\! \sum_{y\in B} \deg(y)(\deg(y)-1) +
        \left( \sum_{x\in A}  \sum_{y \in N(x)} \frac{1}{\deg(x)} \right)
 \left( \sum_{x\in A}  \sum_{y \in N(x)} \deg(x)(\deg(y)-1) \right)
 \end{eqnarray*}
Here the second term is at least
\begin{eqnarray*}
 &\geq&
  \left(\sum_{x\in A}  \sum_{y \in N(x)}
   \frac{\sqrt{\deg(x)(\deg(y)-1)}}{\sqrt{\deg (x)}} \right)^2 =
   \left(\sum_{x\in A}  \sum_{y \in N(x)} \sqrt{\deg(y)-1} \right)^2 \\
 &=&\left(\sum_{y\in B} \deg(y)\sqrt{\deg(y)-1}\right)^2.
  \end{eqnarray*}
Let $F(y_1, y_2, \dots , y_b)$ be a real function
 defined as
$$
  -a\sum_{1\leq i\leq b}(y_i^2-y_i) + \left(\sum_{i}y_i\sqrt{y_i-1}\right)^2,
 $$
where $y_i\geq 2$ and $\sum y_i\geq 2a$.
We obtained that
$aP_3\geq F({\mathbf y})$ where ${\mathbf y}\in R^b$ is the vector
 with coordinates formed by the degrees $\{ \deg(y): y\in B\}$.
We will see that $F$ is convex in this region hence, all $y_i$ can be
 replaced with the average of the degrees, i.e.,  $\sum_{y\in B}
  \deg(y)/b=e/b$. One obtains
$$
  aP_3\geq -ab\frac{e}{b}\left( \frac{e}{b}-1\right)
     +  \left(b\frac{e}{b}\sqrt{\frac{e}{b}-1}\right)^2.
  $$
Rearranging one gets (\ref{eq:P3}).

\noindent{\bf Proof of convexity.}\quad
Let $F_{ij}, F_{ii}$ denote the partial derivatives, $\bf H$ the Hessian of $F$.
Then for $i\neq j$ one has
$$
 F_{ij}=\frac{1}{2}\frac{3y_i-2}{\sqrt{y_i-1}}
 \frac{3y_j-2}{\sqrt{y_j-1}} ,$$
and
\begin{eqnarray*}
 F_{ii}&=&-2a +G({\mathbf y})
 \frac{1}{2}\frac{(3y_i-2)^2}{y_i-1}
  +G({\mathbf y})
 \frac{3y_i-1}{2(y_i-1)\sqrt{y_i-1}}\\
  &\geq& -2a+(G({\mathbf y})-1)
 \frac{1}{2}\frac{(3y_i-2)^2}{y_i-1}
  + \frac{1}{2}\frac{3y_i-2}{\sqrt{y_i-1}}
 \frac{3y_i-2}{\sqrt{y_i-1}}\, ,
  \end{eqnarray*}
where $G({\mathbf y})= \sum y_i\sqrt{y_i-1}$.
Since $(3y-2)^2/2(y-1)\geq 6$ for $y>1$ and $G({\mathbf y})
 \geq \sum {y_i}=e\geq 2a$ we can write $\mathbf H$ as a sum of
 a positive semidefinite matrices, namely $1/2$ times the tensor
 product of the vector $\{\frac{3y_i-2}{\sqrt{y_i-1}} \}$
 with itself,   and a diagonal matrix
 with diagonal entries are exceeding $-2a+ 6(G({\mathbf y})-1)$,
 again a positive definite matrix.
Thus $\mathbf H$ is positive definite and then $F$ is convex in the
 region.
\qed

The above theorem is a slightly improved version of a result of
 Sidorenko~\cite{Sid91} which states
 $P_3\geq e^3/ab- \Delta e$, where $\Delta$ is the maximum degree of $G$.
Concerning general (non-bipartite) graphs Theorem~\ref{W3th} implies
 that $P_3\geq \frac{1}{2}n (d_{\rm ave})^3- \frac{3}{2}n\Delta d_{\rm ave}$.
This inequality may also be deduced from a Moore-type bound, established
 by Alon, Hoory and Linial~\cite{AHL}.

\section{Graphs without $C_6$}

\begin{theorem}\label{thg8bip}
Let $G(A,B)$ be a bipartite graph with $e$ edges and with color classes
 $A$ and $B$, $|A|=a$, $|B|=b$.
Suppose that $G$ has girth eight.
Then for the number edges one has
\begin{equation}\label{eq:g8bip}
  e\leq   (ab)^{2/3} +a +b.
 \end{equation}
 \end{theorem}

\noindent{\it Proof.}\quad
We use induction on the number of vertices if
 there is any isolated vertex, or a vertex of degree 1.
Otherwise, observe, that every pair $x\in A$, $y\in B$ is connected
 by at most one path of length 3.
Thus $P_3\leq ab$.
Comparing this to the lower bound for $P_3$ in
 (\ref{eq:P3}) and rearranging we get the Theorem. \qed
\medskip

D.~de Caen and Sz\'ekely~\cite{dCSz91} showed that
 $e(G)=O((ab)^{2/3})$ assuming $a=O(b^2)$ and $b=O(a^2)$.
Later they showed~\cite{dCSz97} that
 if $G$ has girth eight and every vertex has degree at least
 two, then $e\leq 2^{1/3}(ab)^{2/3}$ and here the
 coefficient $2^{1/3}$ is the best possible by exhibiting a
 graph with $a=2s$, $b=s^2$ and $e=2s^2$.
(Note that this does not contradicts to our result
 (\ref{eq:g8bip}) since here $b=e/2$).

Gy\H ori~\cite{Gyori} observed that in a $C_6$-free graph $G$
 the maximal complete bipartite graphs $K_{\alpha ,\beta}$'s with
 $\alpha,\beta\geq 2$ are edge disjoint  (indeed, these are
 $K_{2,\beta}$'s).
Thus one can remove edges from $G$ such that the resulting graph
 $G_0$ is $C_4$-free and $e(G_0)\geq \frac{1}{2}e(G)$.
Thus Gy\H ori's result combined with Theorem~\ref{thg8bip} gives that
\begin{corollary}\label{th:bipC6}
  If $G$ is a $C_6$-free bipartite graph with
parts of sizes $a$ and $b$ then $e(G)\leq 2(ab)^{2/3}+2a+2b$. \qed
 \end{corollary}

More is true.
In~\cite{C6} it was proved that for such a graph
\begin{equation}\label{bestbip}
  e(G)< 2^{1/3}(ab)^{2/3}+16(a+b)
 \end{equation}
 holds.
Moreover infinitely many examples show that the coefficient $2^{1/3}$ in the
 best possible for large $a$ and $b$ with $b=2a$.

Concerning general (not necessarily bipartite) graphs, it was proved
 by Bondy and Simonovits~\cite{BS} in 1974
 that a graph on $n$ vertices with at least $100kn^{1+1/k}$
 edges contains $C_{2k}$, a cycle of length $2k$.
This was extended into bipartite graphs with parts of sizes
 of $a$ and $b$ by G.~N.~S\'ark\"ozy~\cite{Sarkozy}
 who shoved that such a graph with
 $\max\{ 90k(a+b), 20k(ab)^{1+1/k}\}$ edges contains a $C_{2k}$.
Our Corollary~\ref{th:bipC6} gives these for $C_6$, even a slightly better statement,
 using the following important reduction theorem.

\begin{lemma}\label{th:Erdos} {\rm (Erd\H os~\cite{ErdosBip})}\quad
Let $G$ be an arbitrary graph.
Then there exists a bipartite subgraph $G_0$ with
 $\deg_{G_0}(x)\geq \frac{1}{2}\deg_G(x)$ for all vertices.
Especially, $e(G_0)\geq \frac{1}{2}e(G)$.
 \end{lemma}

\begin{corollary}\label{th:C6}
If $G$ is a $C_6$-free graph on $n$ vertices
 then $e(G)\leq 2^{2/3}n^{4/3}+4n $. \qed
 \end{corollary}

It is known that there are $C_6$-free graphs with at least
 $(\frac{1}{2}+o(1))n^{4/3}$ edges~\cite{Laz}, and
the best known lower and upper bounds can be found in~\cite{C6},
(namely $0.533 n^{4/3}< \ex(n, C_6)< 0.628 n^{4/3}$ for $n> n_0$).
Yuansheng and Rowlinson~\cite{YuangRowliC6}
 determined $\ex(n,C_6)$ and all extremal graphs for $n\le 26$.

\section{Cube-free graphs}

\begin{theorem}[Erd\H os and Simonovits~\cite{ESkocka}]\label{csakQ8th}
Let $Q$ denote the $8$-vertex graph formed by the $12$ edges of a cube.
Then $\ex(n,Q)\leq O(n^{8/5})$.
 \end{theorem}

The original proof of this is rather complicated.
It applies a remarkable regularization process for non-dense
 bipartite graphs.
A somewhat simpler proof was found by Pinchasi and
 Sharir~\cite{PS}, who were interested in certain geometric
 incidence problems, and who also extended it to a bipartite version
\begin{equation}\label{eq:PS}
e(G(A,B))\leq O((ab)^{4/5} +ab^{1/2}+a^{1/2}b).
  \end{equation}
Here we give an even simpler proof which also gives the
 bipartite version, see (\ref{eq:Q8ab}) below.
We only use Theorem \ref{P3th}, Corollary~\ref{th:bipC6} and the power
 mean inequality (\ref{eq:3}), but the main ideas are the same as in~\cite{ESkocka}.

\medskip
\noindent{\bf Proof} of Theorem~\ref{csakQ8th}.\quad
Let $G$ be an $n$-vertex $Q$-free graph.
First, applying Erd\H os' Lemma~\ref{th:Erdos} we choose a
 large bipartite subgraph $G(A,B)$ of $G$, $e(G)\leq 2 e(G(A,B))$.

We say that a hexagon $z_1z_2\dots z_6$ lies {\bf between} the vertices
 $x$ and $y$ if $z_1,z_3,z_5$ are neighbors of $x$ and the other
 vertices of the hexagon are neighbors of $y$, i.e.,
 $z_1,z_3,z_5\in N(x)$ and $z_2,z_4,z_6\in N(y)$ and $\{ x,y \}
 \cap\{ z_1, \dots, z_6\}=\emptyset$.
The crucial observation is that  $x$ and $y$ together with the 6 vertices
 of a hexagon between them contain a cube $Q$.
So there is no such hexagon in a $Q$-free graph.
Thus we can apply the upper bound for the Tur\'an numbers of $C_6$, i.e.,
 Theorem~\ref{th:bipC6} and obtain an upper bound
 for the number of edges $uv$, $u\in N(x)$, $v\in N(y)$.
This gives an upper bound for the number of paths with endvertices $x$ and $y$.
For given $x$ and $y$ we have
\begin{eqnarray*}
  \# xuvy\,\,\, {\rm paths}\,\,
&=& | \{ uv\in E(G(A,B)): u\in N(x)\setminus \{ y\}, v\in N(y)\setminus \{ x\}\}  \\
 &\leq &2|N(x)|^{2/3}|N(y)|^{2/3} + 2|N(x)|+2|N(y)|.
  \end{eqnarray*}
Add this up for every $x\in A$, $y\in B$.
Let $e:=e(G[A,B])$ and use (\ref{eq:3}) with $(r,s)=(1,3/2)$.
We have
\begin{eqnarray*}
  P_3(G(A,B))&\leq&\sum_{x\in A}\sum_{y\in B}
  2\deg(x)^{2/3}\deg(y)^{2/3} + 2\deg(x)+2\deg(y)\\
   &=&2\left( \sum_{x\in A}\deg(x)^{2/3}\right)
 \left( \sum_{y\in B}\deg(y)^{2/3}\right) +2be+2ae\\
 &\leq& 2\times e^{2/3}a^{1/3}\times e^{2/3}b^{1/3}+2(a+b)e.
  \end{eqnarray*}
Comparing this to the lower bound in Theorem~\ref{P3th}
 one obtains that
$$
  (e-a)(e-b)\leq 2e^{1/3}(ab)^{4/3}+2(a+b)ab.
  $$
This implies that
\begin{equation}\label{eq:Q8ab}
  e\leq 2^{3/5}(ab)^{4/5}+2ab^{1/2}+2a^{1/2}b.
  \end{equation}
Using $e(G)\leq 2e$ we obtain
\begin{equation}\label{eq:11}
  e(G)\leq n^{8/5}+(2n)^{3/2}
  \end{equation}
 finishing the proof.
\qed
\medskip

If we use (\ref{bestbip}) instead of  Corollary~\ref{th:bipC6} then the above calculation
 give
 \begin{theorem}
\begin{equation}\label{eq:newPS}
 \ex(a,b,Q)\leq 2^{1/5}(ab)^{4/5} + 9(ab^{1/2}+a^{1/2}b)
  \end{equation}
and
\begin{equation}\label{eq:12}
  \ex(n, Q)\leq 2^{-2/5}n^{8/5}+ 13n^{3/2}.
  \end{equation}
 \end{theorem}

\section{A lower bound on the number of $C_4$'s}

Let $N(G,H)$ denote the number of subgraphs of $G$ isomorphic to $H$.
This function is even more important than the original Tur\'an's problem.
Here we consider only one the simplest cases, $H=C_4$.

It was observed and used many times that for sufficiently large $e$ the graph
 $G$ contains at least $\Omega(e^4/n^4)$ copies of $C_4$.
This result goes back to Erd\H os (1962) and was published, e.g.,
 in Erd\H{o}s and Simonovits~\cite{ESkocka} in an asymptotic form
 ($N(G,C_4)> (1/3)e^4/n^4$ for $n > Cn^{3/2}$).
The following simple form has the advantage that it is valid for
arbitrary $n$ and $e$.

\begin{lemma}[see \cite{FureOzka}]\label{thm-c4}
Let $G$ be a graph with $e$ edges and $n$ vertices. Then
\begin{equation}\label{eq1}
N(G,C_4) \geq  2\frac{e^3(e-n)}{n^4}-\frac{e^2}{2n}\geq
2\frac{e^4}{n^4} - \frac{3}{4}en.
   \end{equation}
\end{lemma}

Allen, Keevash, Sudakov, and Verstra\"ete~\cite{AKSV} gave a bipartite version of
 Lemma~\ref{thm-c4}.
Here we state their result in a slightly stronger form (it is valid for all
 vales of $a,b$ and $e$).
Note that the formula is not symmetric in $A$ and $B$.

\begin{lemma}\label{le:10}
Let $G$ be a bipartite graph with parts $A$ and $B$ of sizes $a$ and $b$ and $e$ edges.
Then the number of $4$-cycles in $G$ is at least
 \begin{equation}\label{eq:AKSV} \frac{e^2(e-b)^2-e(e-b)ba(a-1)}
         {4b^2a(a-1)}.\end{equation}
 \end{lemma}

For completeness we present the proofs of the above Lemmas (below and in the Appendix).
But we will need a slightly stronger and more technical version.

\begin{lemma}\label{le:11}
Let $G$ be a bipartite graph with parts $A$ and $B$ of sizes $a$ and $b$ and $e$ edges.
Let $D(x)$ denote $\sum_{y\in N(x)} (\deg(y)-1)$.
Then the number of $4$-cycles in $G$ is at least
 \begin{equation}\label{eq:bipC4UJ}
  \frac{1}{4(a-1)}\left(\sum_{x\in A} D(x)^2\right) -\frac{1}{4}\left(\sum_{x\in A} D(x)\right).
          \end{equation}\end{lemma}

\medskip
\noindent{\it Proof.}\quad
We have
\begin{eqnarray*}
 N(G, C_4)&=&\sum_{\{x,x'\}\subset  A}  {d(x,x')\choose 2}=
   \frac{1}{2}\sum_{x\in A}\left(\sum_{x'\in A\setminus x}  {d(x,x')\choose 2}\right)\\
   &=& \frac{1}{2}\sum_{x\in A} (a-1)  { \sum_{x'\in A\setminus x} d(x,x')/(a-1)\choose 2}\\
 &=&  \frac{a-1}{2}\sum_{x\in A}  { D(x)/(a-1)\choose 2}\\
  &=&  \frac{1}{4(a-1)}\left(\sum_{x\in A} D(x)^2\right) -\frac{1}{4}\left(\sum_{x\in A} D(x)\right).\qed
  \end{eqnarray*}

Note that Lemma~\ref{le:11} easily implies Lemma~\ref{le:10}.
Indeed, observe that for $e(e-b)<  ba(a-1)$ the right hand side of (\ref{eq:AKSV}) is negative,
 so we may suppose that  $(e^2/b) - e\geq  a(a-1)$.
Use Chauchy-Schwartz for $\sum_{x\in A} D(x)$. We obtain
$$ \sum_{x\in A} D(x)= \sum_{x\in A} \left(\sum_{y\in B, xy\in E(G)} (\deg(y)-1)\right)=\sum_{y\in B}
  \deg(y)^2 -\sum_{y\in B} \deg(y) \geq \frac{e^2}{b}-e.
  $$
Use Chauchy-Schwartz again for $\sum D(x)^2$. We have
$$ \sum_{x\in A} D(x)^2\geq\frac{1}{a} \left( \sum_{x\in A} D(x) \right)^2.
  $$
Now Lemma~\ref{le:11} gives that
  $N(G,C_4)\geq (N^2/4a(a-1))-(N/4)$ for $N:=\sum_{x\in A} D(x)$.
Since $N\geq (e^2/b) - e\geq a(a-1)$ the polynomial $p(N):=N^2/a(a-1)-N$ is increasing and
 we get $N(G,C_4)\geq p(N)\geq p(e(e-b)/b))$. \qed

\section{Cubes with a diagonal}

\begin{theorem}[Erd\H os and Simonovits~\cite{ESkocka}]\label{Q8+th}
Let $Q^+$ denote the $8$-vertex graph formed by the $12$ edges of a cube
 with a long diagonal.
Then $\ex(n,Q^+)\leq O(n^{8/5})$.
 \end{theorem}

Here we give a simpler proof which also gives a stronger
 bipartite version.

\begin{theorem}\label{th:Q8+Bip}
${}$ \quad $\ex(a,b,Q^+)= 2^{3/5} (ab)^{4/5}+ O(ab^{1/2}+a^{1/2}b )$.
 \end{theorem}
Using again Erd\H os' Lemma~\ref{th:Erdos} and $a+b=n$ we get
\begin{equation}\label{eq:jav}
  \ex(n,Q^+)\leq  n^{8/5}+ O(n^{3/2}).
  \end{equation}

\medskip
\noindent{\it Proof of Theorem~\ref{th:Q8+Bip}}.\quad
Let $G$ be an $n$-vertex $Q^+$-free bipartite graph with classes $A$ and $B$.
The main idea is the same as in~\cite{ESkocka} and in the proof of Theorem~\ref{csakQ8th}.
The crucial observation is that an edge $xy\in E(G)$  together with
 the 6 vertices of a hexagon between them form a $Q^+$.
So there is no such hexagon in a $Q^+$-free graph between the neighborhoods
 of two connected vertices.
Thus we can apply the upper bound for the Tur\'an numbers of $C_6$, i.e.,
 Theorem~\ref{th:bipC6} and obtain an upper bound
 for the number of edges $x'y'$, $y'\in N(x)$, $x'\in N(y)$ for $xy\in E(G)$.
This gives an upper bound for the number of fourcycles containing the edge $xy$.
\begin{eqnarray*}
  \# xy'x'y\,\,\, {\rm four\, \,  cycles}\,\,
&=& | \{ x'y'\in E(G(A,B)): y'\in N(x)\setminus \{ y\}, x'\in N(y)\setminus \{ x\}\}  \\
 &\leq &2(|N(x)|-1)^{2/3}(|N(y)|-1)^{2/3} + 2|N(x)|+2|N(y)|-2.
  \end{eqnarray*}
Add this up for every $x\in A$, $y\in B$, $xy\in E(G)$ and apply  (\ref{eq:3})
 for $\sum_{xy\in E}(\deg(y)-1)^{2/3}$ with $(r,s)=(1,3/2)$ for every $x$. We obtain
\begin{eqnarray*}
  4N(G,C_4) &\leq& \sum_{x\in A}\sum_{y\in N(x)} 2(\deg(x)-1)^{2/3}(\deg(y)-1)^{2/3} + 2\deg(x)+2\deg(y)-2\\
    &=&2\left(   \sum_{x\in A}  (\deg(x)-1)^{2/3}\deg(x)^{1/3}D(x)^{2/3}  \right)
    +2\left( \sum_{y\in B}  D(y)  \right)+2\left( \sum_{x\in A} D(x)\right).
    \end{eqnarray*}
Apply H\"older inequality with $1/p=2/3$ and $1/q=1/3$ in the first term. We obtain
 it is at most
\begin{equation}\label{eq:17}
\leq 2\left(   \sum_{x\in A}  (\deg(x)-1)\deg(x)^{1/2}\right)^{2/3}\left(   \sum_{x\in A}D(x)^{2}  \right)^{1/3}.
    \end{equation}

From now on, to save time and energy, and to better emphasize the main steps of calculation
  we only sketch the proof.
Compare the obvious leading terms in the lower and upper bounds (\ref{eq:bipC4UJ}) and (\ref{eq:17})
 for $N(G,C_4)$, we have
$$
\frac{1}{a-1}\left(   \sum_{x\in A}D(x)^{2}  \right)   \ll 4N(G,C_4)\ll 2\left(   \sum_{x\in A}  (\deg(x)-1)\deg(x)^{1/2}\right)^{2/3}\left(   \sum_{x\in A}D(x)^{2}  \right)^{1/3}
    $$
yielding
\begin{equation}\label{eq:18}
 \left(   \sum_{x\in A}D(x)^{2}  \right)   \ll(2(a-1))^{3/2}\left(   \sum_{x\in A}  (\deg(x)-1)\deg(x)^{1/2}\right).
    \end{equation}
On the left hand side we can use Cauchy-Schwartz  and
  on the right hand side we apply  (\ref{eq:3}) with $(r,s)=(3/2,2)$. We obtain
$$
 \frac{1}{a}\left(   \sum_{x\in A}D(x)  \right)^2
   \ll(2(a-1))^{3/2} a^{1/4}\left(   \sum_{x\in A}  (\deg(x)-1)d(x)\right)^{3/4}.
    $$
Rearranging we have
  \begin{equation}\label{eq:19}
  \left(   \sum_{x\in A}D(x)  \right)^2
   \ll 2^{3/2} a^{11/4}\left(   \sum_{y\in B} D(y)\right)^{3/4}.
    \end{equation}

Exchange the role of $A$ and $B$, we get
$$
  \left(   \sum_{x\in A}D(x)  \right)^2
   \ll 2^{3/2} b^{11/4}\left(   \sum_{x\in A} D(x)\right)^{3/4}.
    $$
Multiply the above two inequalities, take 4th power, we get
$$
  2^{12} (ab)^{11}\gg \left(   \sum_{x\in A}D(x)  \right)^5\left(   \sum_{y\in B}D(y)  \right)^5
   \ge \left( \frac{e^2}{a}-e\right)^5\left( \frac{e^2}{b}-e\right)^5
$$
 leading to $2^{12}(ab)^{16} \gg e^{20}$. \qed

\section{Appendix}\label{S:App}

1.\quad  \noindent{\it A direct proof of the Mulholland-Smith inequality (\ref{eq:5}) concerning
 the number of 3-walks using only high school algebra.}

Considering the middle edge of the 3-walks one obtains that
$$
 W_3=\sum_{x\in V}  \sum_{y \in N(x)} \deg(x)\deg(y) .
 $$
Here we have $2e=nd_{\rm ave}$ terms.
Our aim is to separate the variables in the products $\deg(x)\deg(y)$
 so next we use first that the $2e$-dimensional Quadratic inequality
(Quadratic mean is greater or equal than the Arithmetic mean),
 second we use (for 2 variables) that the Arithmetic mean is
 greater or equal than the Harmonic mean, then third time we use
 again (this time for $2e$ variables) that Arithmetic $\geq$
 Harmonic.
One obtains that
\begin{eqnarray*} 
\sqrt{ \frac{W_3}{nd_{\rm ave}}} &=& \sqrt{\frac{\sum_{x\in V}  \sum_{y \in
      N(x)} \deg(x)\deg(y)}{2e}}\\
  &\geq& \frac{\sum_{x\in V}  \sum_{y \in N(x)} \sqrt{\deg(x)\deg(y)}}{2e} \\
  &\geq& \frac{1}{2e}\left(  \sum_{x\in V}  \sum_{y \in N(x)}
    \frac{2}{\frac{1}{\deg(x)}+ \frac{1}{\deg(y)}} \right) \\
   &\geq&2e\left( \sum_{x\in V}  \sum_{y \in N(x)}
 \frac{ \frac{1}{\deg(x)}+ \frac{1}{\deg(y)}}{2} \right)^{-1}
 = \frac{2e}{\sum_{x\in V} 1} =\frac{2e}{n}=d_{\rm ave}.
  \end{eqnarray*}%
\qed

\medskip
\noindent
2.\quad  {\it Proof of Lemma~\ref{thm-c4} concerning the number of $C_4$'s.}

 Denote the number of $x,y$-paths of length two by $d(x,y)$.
We have
\begin{equation}\label{eq2}
{\overline{\overline{d}}}:={n \choose 2}^{-1}\sum_{x,y\in V(G)}d(x,y) = {n \choose 2}^{-1} \sum_{x\in
V(G)} {\deg(x) \choose 2} \geq {n \choose 2}^{-1} n {2e/n \choose 2}.
  \end{equation}
Therefore,
${\overline{\overline{d}}} \geq \frac{2e(2e-n)}{n^2(n-1)}$.
Moreover
\begin{equation}\label{eq3}
 N(G,C_4) = \frac{1}{2} \sum_{x,y\in V(G)} {d(x,y) \choose 2} \geq \frac{1}{2} {n
\choose 2}{{\overline{\overline{d}}} \choose 2}.
  \end{equation}
We may suppose that the middle term in (\ref{eq1})
 is positive, which implies that $\frac{2e(2e-n)}{n^2(n-1)}\geq 1/2$.
The paraboloid ${x\choose 2}$ is increasing for $x\geq 1/2$.
So we may substitute the lower bound of ${\overline{\overline{d}}}$ from (\ref{eq2})
 into (\ref{eq3}) and a little algebra gives (\ref{eq1}).
\qed

\end{document}